\documentclass{article}\usepackage{amsmath,amssymb}
\newtheorem{definition}{Definition}[section]\newtheorem{theorem}[definition]{Theorem}
\newtheorem{corollary}[definition]{Corollary}\newtheorem{lemma}[definition]{Lemma}\newtheorem{remark}[definition]{Remark}\newtheorem{example}{Example}\newtheorem{proof}{Proof}
\begin{document}

\begin{center}
{\large\bf Right-covariant differential calculus on ${\cal F}({\mathbb C}_q^{2|1})$}

\vspace*{.3cm}
Salih Celik

\vspace*{.3cm}
Department of Mathematics, Yildiz Technical University, DAVUTPASA-Esenler,
Istanbul, 34210 TURKEY.
\end{center}

\begin{abstract}
We define a new ${\mathbb Z}_2$-graded quantum (2+1)-space and show that the extended ${\mathbb Z}_2$-graded algebra of polynomials on this
${\mathbb Z}_2$-graded quantum space, denoted by ${\cal F}({\mathbb C}_q^{2\vert1})$, is a ${\mathbb Z}_2$-graded Hopf algebra. We construct a right-covariant differential calculus on ${\cal F}({\mathbb C}_q^{2\vert1})$ and define a ${\mathbb Z}_2$-graded quantum Weyl algebra and mention a few algebraic properties of this algebra. Finally, we explicitly construct the dual ${\mathbb Z}_2$-graded Hopf algebra of ${\cal F}({\mathbb C}_q^{2\vert1})$.
\end{abstract}


\vspace*{.3cm}

\section{Introduction}
\label{sec1}

The noncommutative differential geometry of quantum groups was first introduced by Woronowicz in \cite{20} and \cite{21}. Differential calculi on certain classes of quantum homogeneous spaces are described in \cite{2,12,16}. The quantum plane and superplane are the simplest samples of a noncommutative spaces. Following Woronowicz's approach, noncommutative differential calculi on some lower dimensional (super)spaces are presented in \cite{3,4,5,6,7,8}.

The other approach, initiated by Wess and Zumino \cite{19}, followed Manin's emphasis on the quantum spaces \cite{14} as the primary objects. Differential forms are defined in terms of noncommuting coordinates, and the differential and algebraic properties of quantum groups acting on these spaces are obtained from the properties of the spaces. The natural extension of their scheme to ${\mathbb Z}_2$-graded space \cite{15} was introduced in \cite{9,11,13,17}.

In this paper, we introduce ${\mathbb Z}_2$-graded Hopf algebra structure on the extended function algebra of the ${\mathbb Z}_2$-graded quantum space ${\mathbb C}_q^{2\vert1}$ and investigate its noncommutative geometry. In Section 3, we define a ${\mathbb Z}_2$-graded quantum (2+1)-space ${\mathbb C}_q^{2\vert1}$ and show that the ${\mathbb Z}_2$-graded algebra obtained by extending the algebra of functions on this ${\mathbb Z}_2$-graded space has a ${\mathbb Z}_2$-graded Hopf algebra structure. In Section 4, we build a right-covariant differential calculus on ${\cal F}({\mathbb C}_q^{2\vert1})$ and we describe a ${\mathbb Z}_2$-graded quantum Weyl algebra and talk about some of its properties. In the last section, we construct a dual ${\mathbb Z}_2$-graded algebra for ${\cal F}({\mathbb C}_q^{2\vert1})$ within the framework of Sudbery's approach and present the commutation relations of the generators and Hopf structure of this dual ${\mathbb Z}_2$-graded algebra.

\section{Preliminaries}
\label{Sec:2}

In this section, we will briefly talk about some of known basic concepts as far as we need. Throughout the paper, we will fix a base field ${\mathbb C}$, the set of complex numbers. We write ${\mathbb Z}_2 = {\mathbb Z}/2{\mathbb Z}$ and denote the elements of ${\mathbb Z}_2$ by $0,1$.

\subsection{${\mathbb Z}_2$-graded algebras}
\label{Subsec:2.1}

A supervector space $V$ over ${\mathbb C}$ is a ${\mathbb Z}_2$-graded vector space over ${\mathbb C}$ and we write $V =V_0\oplus V_1$, where $V_0$ and $V_1$ are even and odd subspaces of $V$, respectively. The elements of $V_0$ and $V_1$ are called even and odd, respectively. The elements of $V_0\cup V_1$ will be called homogeneous. For a homogeneous element $v$ we write $p(v)$ for the parity or degree; if $v\in V_0$ (resp. $V_1$) we have $p(v)=0$ (resp. $1$).

A superalgebra or ${\mathbb Z}_2$-graded algebra ${\cal A}$ over ${\mathbb C}$ is a supervector space over ${\mathbb C}$ with a map
${\cal A}\times{\cal A}\longrightarrow {\cal A}$ such that ${\cal A}_i \cdot{\cal A}_j \subset {\cal A}_{i+j}$ for $i,j=0,1$. If ${\cal A}$ and ${\cal B}$ are two ${\mathbb Z}_2$-graded algebras then the tensor product ${\cal A}\otimes{\cal B}$ exists and ${\mathbb Z}_2$-graded with
\begin{equation*}
({\cal A}\otimes {\cal B})_i = \bigoplus_{j+k=i} {\cal A}_j\otimes {\cal B}_k, \qquad i\in{\mathbb Z}_2.
\end{equation*}
The following definition gives the product rule for tensor product of ${\mathbb Z}_2$-graded algebras.

\begin{definition} \label{def2.1}
Let ${\cal A}$ and ${\cal B}$ be two ${\mathbb Z}_2$-graded algebras and $a_i$'s and $b_i$'s are homogeneous elements in ${\cal A}$ and ${\cal B}$, respectively. Then, the tensor product rule of the algebras ${\cal A}$ and ${\cal B}$ is defined by
\begin{equation}\label{2.1}
(a_1\otimes a_2)(b_1\otimes b_2)=(-1)^{p(a_2)p(b_1)}(a_1b_1\otimes a_2b_2).
\end{equation}
\end{definition}

\subsection{${\mathbb Z}_2$-graded Hopf algebras}
\label{Subsec:2.2}

The elementary properties of a ${\mathbb Z}_2$-graded Hopf algebra are similar to the corresponding properties of ordinary Hopf algebras.

\begin{definition} \label{def2.2}
A ${\mathbb Z}_2$-graded Hopf algebra (or Hopf superalgebra) is a ${\mathbb Z}_2$-graded vector space ${\cal H}$ over ${\mathbb C}$ with two algebra homomorphisms
$\Delta: {\cal H}\to {\cal H}\otimes {\cal H}$, called the {\it coproduct}, $\epsilon: {\cal H}\to {\mathbb C}$, called the {\it counit} and an algebra antihomomorphism
$S:{\cal H}\to {\cal H}$, called the {\it antipode}, such that
\begin{align} \label{2.2}
(\Delta \otimes {\rm id}) \circ \Delta &= ({\rm id} \otimes \Delta) \circ \Delta, \nonumber\\
m\circ(\epsilon\otimes{\rm id})\circ\Delta &= {\rm id} = m\circ({\rm id}\otimes\epsilon)\circ\Delta, \\
m\circ(S\otimes{\rm id})\circ\Delta &= \eta\circ\epsilon = m\circ({\rm id}\otimes S)\circ\Delta, \nonumber
\end{align}
and $\Delta({\bf 1}) ={\bf 1}\otimes{\bf 1}$, $\epsilon({\bf 1})=1$, $S({\bf 1})={\bf 1}$, where $m$ is the multiplication map, ${\rm id}$ is the identity map and
$\eta:{\mathbb C}\longrightarrow {\cal H}$.
\end{definition}

\begin{remark} \label{rem2.1}
An element of a ${\mathbb Z}_2$-graded Hopf algebra ${\cal H}$ is expressed as a product on the generators and its antipode $S$ is calculated with the property
\begin{equation}\label{2.3}
S(ab) = (-1)^{p(a)p(b)} S(b)S(a), \qquad \forall a,b \in {\cal H}.
\end{equation}
\end{remark}

In many cases, a ${\mathbb Z}_2$-graded Hopf algebra coacts on a ${\mathbb Z}_2$-graded algebra such that the coaction respects the ${\mathbb Z}_2$-graded algebra structure. This situation is explained by the following concept.

\begin{definition} \label{def2.3}
Let ${\cal X}$ be a ${\mathbb Z}_2$-graded algebra and ${\cal H}$ a ${\mathbb Z}_2$-graded Hopf algebra. Then the ${\mathbb Z}_2$-graded algebra ${\cal X}$ is called a left
${\cal H}$-comodule algebra if there exists an algebra homomorphism $\delta_L:{\cal X}\longrightarrow {\cal H} \otimes {\cal X}$ such that
\begin{equation*}
({\rm id}\otimes\delta_L)\circ\delta_L = (\Delta\otimes{\rm id})\circ\delta_L \quad {\rm and} \quad (\epsilon\otimes {\rm id})\circ\delta_L = {\rm id}.
\end{equation*}
Also, $\delta_L$ is called the left coaction of ${\cal H}$ on ${\cal X}$.
\end{definition}

\begin{definition} \label{def2.4}
Let ${\cal X}$ be a ${\mathbb Z}_2$-graded algebra and ${\cal H}$ a ${\mathbb Z}_2$-graded Hopf algebra. Then the ${\mathbb Z}_2$-graded algebra ${\cal X}$ is called a right
${\cal H}$-comodule algebra if there exists an algebra homomorphism $\delta_R:{\cal X}\longrightarrow {\cal X}\otimes {\cal H}$ such that
\begin{equation*}
(\delta_R\otimes{\rm id})\circ\delta_R = ({\rm id}\otimes\Delta)\circ\delta_R \quad {\rm and} \quad ({\rm id}\otimes\epsilon)\circ\delta_R = {\rm id}.
\end{equation*}
Also, $\delta_R$ is called the right coaction of ${\cal H}$ on ${\cal X}$.
\end{definition}

\begin{remark} \label{rem2.2}
A quantum space ${\cal X}$ for ${\cal H}$ is called quantum homogeneous space if there exists an embedding $i_{\cal X}\to {\cal H}$ such that
$\delta_R = \Delta\circ i_{\cal X}$, i.e. ${\cal X}$ can be identified with a subalgebra of ${\cal H}$, and its coaction is obtained then by restricting the coproduct.
\end{remark}

\section{A ${\mathbb Z}_2$-graded Hopf algebra ${\cal F}({\mathbb C}_q^{2\vert1})$}
\label{Sec:3}

The elements of the ${\mathbb Z}_2$-graded (2+1)-space ${\mathbb C}_q^{2\vert1}$ are ${\mathbb Z}_2$-graded vectors (or supervectors) generated by two even and an odd components. In this section, we will introduce a ${\mathbb Z}_2$-graded Hopf algebra structure on the extended function algebra of the ${\mathbb Z}_2$-graded quantum space defined below.

\subsection{The algebra of polynomials on ${\mathbb C}_q^{2\vert1}$}
\label{Subsec:3.1}

Let ${\mathbb C}\langle x_1,x_2,\theta \rangle$ be a free ${\mathbb Z}_2$-graded algebra with unit generated by $x_i$ and $\theta$, where $p(x_i)=0$ for $i=1,2$, and $p(\theta)=1$.

\begin{definition} \label{def3.1}
Let $I_q$ be the two-sided ideal of the algebra ${\mathbb C}\langle x_1,x_2,\theta \rangle$ which generated by the elements $x_1x_2-x_2x_1$, $x_i\theta-q\theta x_i$ and $\theta^2$ for $i=1,2$, where $q\in{\mathbb C}-\{0\}$. The ${\mathbb Z}_2$-graded, associative, unital algebra
$${\cal O}({\mathbb C}_q^{2\vert1}) = {\mathbb C}\langle x_1,x_2,\theta \rangle/I_q$$
is the algebra of polynomials on the ${\mathbb Z}_2$-graded quantum space ${\mathbb C}_q^{2\vert1}$.
\end{definition}
This associative algebra over the complex number is known as the algebra of polynomials over ${\mathbb Z}_2$-graded quantum (2+1)-space. In accordance with
Definition \ref{def3.1}, if $(x_1,x_2,\theta)^t$ belongs to ${\mathbb C}_q^{2\vert1}$ then we have
\begin{equation} \label{3.1}
x_1x_2 = x_2x_1, \quad x_i\theta = q\theta x_i, \quad \theta^2 =0, \qquad (i=1,2)
\end{equation}
where $q$ is a non-zero complex number.

If we consider the generators of the algebra ${\cal O}({\mathbb C}_q^{2\vert1})$ as linear functionals, we can find many 3x3-matrix representations of these generators that preserve the relations (\ref{3.1}):

\begin{example}
(a) It can be seen that the matrices
\begin{equation*}
x_1:=\begin{pmatrix} q & 0 & 0 \\ 0 & q & 0 \\ 0 & 0 & 1 \end{pmatrix}, \quad
x_2:=\begin{pmatrix} 1 & q^2-1 & 0 \\ 0 & 1 & 0 \\ 0 & 0 & q^{-1} \end{pmatrix}, \quad
\theta:=\begin{pmatrix} 0 & 0 & q^2-1 \\ 0 & 0 & 0 \\ 0 & 0 & 0 \end{pmatrix}
\end{equation*}
corresponding to the coordinate functions satisfy relations $(\ref{3.1})$.

\noindent(b) It can be seen that the matrices
\begin{equation*}
x_1:=\begin{pmatrix} 1 & 0 & 0 \\ 0 & 1 & 0 \\ 0 & 0 & q \end{pmatrix}, \quad
x_2:=\begin{pmatrix} 1 & 0 & 0 \\ q^{-2}-1 & 1 & 0 \\ 0 & 0 & q \end{pmatrix}, \quad
\theta:=\begin{pmatrix} 0 & 0 & 0 \\ 0 & 0 & 0 \\ 1-q^{-2} & 0 & 0 \end{pmatrix}
\end{equation*}
corresponding to the coordinate functions satisfy relations $(\ref{3.1})$.
\end{example}

\begin{definition} \label{def3.2}
Let $\Lambda({\mathbb C}_q^{2\vert1})$ be the algebra with the generators $\varphi_1$, $\varphi_2$ and $z$ satisfying the relations
\begin{equation} \label{3.2}
\varphi_i\varphi_j = - q^{2(i-j)} \varphi_j\varphi_i, \quad \varphi_i z = q^{-1} z \varphi_i, \qquad (i,j=1,2)
\end{equation}
where the coordinates $\varphi_i$ are of degree 1, the coordinate $z$ is of degree 0. We call $\Lambda({\mathbb C}_q^{2\vert1})$ exterior algebra of the
${\mathbb Z}_2$-graded space ${\mathbb C}_q^{2\vert1}$.
\end{definition}

\subsection{${\mathbb Z}_2$-graded Hopf algebra ${\cal F}({\mathbb C}_q^{2\vert1})$}
\label{Subsec:3.2}

Let us suppose that the inverses $x^{-1}_i$ of $x_i$ exist such that $x_ix^{-1}_i = {\bf 1} =  x^{-1}_i x_i$ for $i=1,2$. Then, we define the extended ${\mathbb Z}_2$-graded quantum algebra, denoted by ${\cal F}({\mathbb C}_q^{2\vert1})$, to be the ${\mathbb Z}_2$-graded algebra containing ${\cal O}({\mathbb C}_q^{2\vert1})$ with the inverses $x^{-1}_i$ and the unit {\bf 1}.

\begin{theorem} \label{thm3.1}
The ${\mathbb Z}_2$-graded algebra ${\cal F}({\mathbb C}_q^{2\vert1})$ is a ${\mathbb Z}_2$-graded Hopf algebra. The definitions related to a co-maps on
${\cal F}({\mathbb C}_q^{2\vert1})$ are as follows:

\noindent(i) The coproduct $\Delta: {\cal F}({\mathbb C}_q^{2\vert1}) \longrightarrow {\cal F}({\mathbb C}_q^{2\vert1}) \otimes {\cal F}({\mathbb C}_q^{2\vert1})$ is given by
\begin{equation} \label{3.3}
\Delta(x_i) = x_i \otimes x_i, \quad \Delta(\theta) = \theta \otimes x^{-1}_1x_2 + {\bf 1} \otimes \theta, \qquad (i=1,2).
\end{equation}

\noindent(ii) The counit $\epsilon: {\cal F}({\mathbb C}_q^{2\vert1}) \longrightarrow {\mathbb C}$ is defined by
\begin{equation*}
\epsilon(x_i)=1, \quad \epsilon(\theta)=0, \qquad (i=1,2).
\end{equation*}

\noindent(iii) The ${\mathbb Z}_2$-graded algebra ${\cal F}({\mathbb C}_q^{2\vert1})$ admits a coinverse $S:{\cal F}({\mathbb C}_q^{2\vert1}) \longrightarrow F({\mathbb C}_{q^{-1}}^{2\vert1})$ defined by
\begin{equation*}
S(x_i) = x^{-1}_i, \quad S(\theta) = - q^{-1} x_1\theta x^{-1}_2, \qquad (i=1,2).
\end{equation*}
\end{theorem}

\begin{proof}
It is necessary to show that the coproduct, counit and coinverse described above define a ${\mathbb Z}_2$-graded Hopf algebra. For this, it is sufficient to check the conditions given in Definition \ref{def2.2} on each generator separately. These are carried out  by making direct calculations.
\end{proof}

The set $\{x^k_1x_2^l\theta^m: \, k,l\in{\mathbb N}_0, m=0,1\}$ forms a vector space basis of ${\cal O}({\mathbb C}_q^{2\vert1})$. Since $\Delta$ is a ${\mathbb Z}_2$-graded algebra homomorphism, we can easily calculate the action of $\Delta$ on all sums of products of generators can be computed.

\subsection{$*$-Structure on the algebra ${\cal F}({\mathbb C}_q^{2\vert1})$}
\label{Sec:3.3}

Here, we define a ${\mathbb Z}_2$-graded involution (or superinvolution) on the algebra ${\cal F}({\mathbb C}_q^{2\vert1})$ and show that ${\cal F}({\mathbb C}_q^{2\vert1})$ is a
${\mathbb Z}_2$-graded $*$-bialgebra.

\begin{definition} \label{def3.3}
Let ${\cal A}$ be an associative ${\mathbb Z}_2$-graded algebra. A map $*: {\cal A}\longrightarrow{\cal A}$ is called a ${\mathbb Z}_2$-graded involution if
\begin{equation*}
(ab)^*=(-1)^{p(a)p(b)} b^* a^*, \qquad (a^*)^*=a
\end{equation*}
for any elements $a,b\in{\cal A}$. A homogeneous element $a\in{\cal A}$ is called supersymmetric if $a^*=(-1)^{p(a)} a$. The pair $({\cal A},*)$ is called a ${\mathbb Z}_2$-graded $*$-algebra.
\end{definition}

\begin{theorem} \label{thm3.2}
If $q\ne\pm1$ is a complex number of modulus one, then the algebra ${\cal O}({\mathbb C}_q^{2\vert1})$ equipped with the ${\mathbb Z}_2$-graded involution determined by
\begin{equation} \label{3.4}
x^*_i := x_i, \quad \theta^* := -\theta \qquad (i=1,2)
\end{equation}
becomes a ${\mathbb Z}_2$-graded $*$-algebra ${\cal O}({\mathbb R}_q^{2\vert1})$.
\end{theorem}

\begin{proof}
If $q\ne\pm1$, using the relations in (\ref{3.1}), we write
\begin{equation*}
0 = (x_i\theta - q\theta x_i)^* = -\theta x_i + \bar{q} x_i\theta = (q\bar{q}-1) \theta x_i
\end{equation*}
which gives $q\bar{q}=1$.
\end{proof}

If ${\cal A}$ is a $*$-algebra then the involution of ${\cal A}\otimes{\cal A}$ is defined by $(a\otimes b)^*=a^*\otimes b^*$ for all elements
$a,b \in {\cal A}$.

\begin{definition} \label{def3.4}
A ${\mathbb Z}_2$-graded bialgebra ${\cal A}$ is called a ${\mathbb Z}_2$-graded $*$-bialgebra if ${\cal A}$ is equipped with an operation $*$ such that
\begin{equation} \label{3.5}
\Delta(a^*)=[\Delta(a)]^*, \quad \epsilon(a^*)=\overline{\epsilon(a)}, \quad a\in{\cal A}. 
\end{equation}
A ${\mathbb Z}_2$-graded Hopf algebra which is a $*$-bialgebra is called a ${\mathbb Z}_2$-graded Hopf $*$-algebra.
\end{definition}

\begin{theorem} \label{thm3.3}
If $\bar{q}=q^{-1}$ then the ${\mathbb Z}_2$-graded algebra ${\cal F}({\mathbb C}_q^{2\vert1})$ equipped with the involution given by $(\ref{3.4})$
becomes a ${\mathbb Z}_2$-graded $*$-bialgebra ${\cal F}({\mathbb R}_q^{2\vert1})$.
\end{theorem}

\begin{proof} In accordance with Definition \ref{def3.4}, one should show that the equalities in (\ref{3.5}) are fulfilled. These can be shown by direct calculation.
\end{proof}

\section{Right-covariant differential calculus on ${\cal F}({\mathbb C}_q^{2\vert1})$}
\label{Sec:4}

It is well known that functions commute with differentials or 1-forms in usual algebraic geometry. That is, the multiplication of functions from the left or from the right by differentials is the same. However, this case cannot be true when the algebra of functions is noncommutative. A differential calculus is said to be noncommutative if the left and right multiplication of functions by 1-forms do not coincide.

In this section, we will give right-covariant first order ${\mathbb Z}_2$-graded differential calculus on the ${\mathbb Z}_2$-graded Hopf algebra
${\cal F}({\mathbb C}_q^{2\vert1})$ generated by three generators and quadratic relations. This calculus involves functions on the ${\mathbb Z}_2$-graded space
${\mathbb C}_q^{2\vert1}$ and differentials. Thus a linear operator {\sf d} which acts on the generators of the algebra ${\cal F}({\mathbb C}_q^{2\vert1})$ must be defined. For the definition, it is sufficient to define an action of {\sf d} on the generators and on their products as follows: for example, the operator {\sf d} is applied to $x_1$ or $x_2$ produces a 1-form of degree 1, by the definition. Similarly, application of {\sf d} to $\theta$ produces 1-form whose ${\mathbb Z}_2$-grade is 0. In order to set up a differential scheme, we will use the right-covariance.

\subsection{First order differential calculi on ${\cal F}({\mathbb C}_q^{2\vert1})$}
\label{Subsec:4.1}

Using the general theory formulated by Woronowicz in \cite{21} (see, also \cite{2}), we will recall some basic notions about differential calculus on ${\mathbb Z}_2$-graded Hopf algebra. Let us begin with the definition of a first order ${\mathbb Z}_2$-graded differential calculus on a general algebra ${\cal A}$ and specify later to the case of a
${\mathbb Z}_2$-graded Hopf algebra.

\begin{definition} \label{def4.1} 
Let ${\cal A}$ be an associative ${\mathbb Z}_2$-graded algebra and $\Gamma$ an ${\cal A}$-bimodule. Then the pair $(\Gamma,{\sf d})$ is called a first order
${\mathbb Z}_2$-graded differential calculus over ${\cal A}$ if ${\sf d}:{\cal A}\longrightarrow\Gamma$ is a linear map satisfying the
${\mathbb Z}_2$-graded Leibniz rule
\begin{align*}
{\sf d}(a\cdot b) &= ({\sf d}a)\cdot b + (-1)^{p(a)} \, a\cdot {\sf d}b, \quad \forall a,b\in{\cal A},
\end{align*}
such that $\Gamma={\rm Lin}\{a\cdot {\sf d}b: \, a,b \in {\cal A}\}$. The elements of $\Gamma$ is so called 1-forms. Any element $w\in\Gamma$ can be written
$w=\sum_k a_k\cdot {\sf d}b_k$, where $a_k,b_k\in {\cal A}$.
\end{definition}

\begin{definition} \label{def4.2} 
Let $(\Gamma,{\sf d})$ be a first order ${\mathbb Z}_2$-graded differential calculus over ${\mathbb Z}_2$-graded $*$-algebra ${\cal A}$. Then $(\Gamma,{\sf d})$ is called a $*$-calculus if there exists an involution $*:\Gamma\longrightarrow\Gamma$ such that
\begin{equation*}
(a\cdot{\sf d}b)^*=(-1)^{p(a)p({\sf d}b)}\, {\sf d}(b^*)\cdot a^*
\end{equation*}
for $a,b\in{\cal A}$.
\end{definition}

The right-covariant differential calculus is playing a critical role in the analysis of differential structures on ${\mathbb Z}_2$-graded Hopf algebras. Let ${\cal H}$ be a
${\mathbb Z}_2$-graded Hopf algebra with unit element {\bf 1} and let $(\Gamma,{\sf d})$ be a first order ${\mathbb Z}_2$-graded differential calculus over ${\cal H}$.

\begin{definition} \label{def4.3} 
The differential calculus $(\Gamma,{\sf d})$ is called right-covariant if there is a linear map $\Delta_R:\Gamma\longrightarrow \Gamma\otimes {\cal H}$, which is  called right coaction, such that
\begin{equation} \label{4.1}
\Delta_R(a \cdot {\sf d}b) = \Delta(a)({\sf d}\otimes{\rm id})\Delta(b), \quad \forall a,b\in {\cal H}.
\end{equation}
\end{definition}

To obtain $q$-commutation relations between the generators of ${\cal F}({\mathbb C}_q^{2\vert1})$ and their first order differentials, we combine $x_1$, $x_2$, $\theta$ and their differentials ${\sf d}x_1$, ${\sf d}x_2$, ${\sf d}\theta$ which are considered as elements of a space
$\Gamma^1 := \Gamma^1_R({\cal F}({\mathbb C}_q^{2\vert1}))$ of 1-forms. We allow multiplication of the first order differentials in $\Gamma^1$ from the left and from the right by the elements of ${\cal F}({\mathbb C}_q^{2\vert1})$, so that by the definition of the multiplication the resulting 1-form belongs to $\Gamma^1$. This means that $\Gamma^1$ is an
${\cal F}({\mathbb C}_q^{2\vert1})$-bimodule.

Let ${\cal F}({\mathbb C}_q^{2\vert1})$-bimodule $\Gamma^1$ be generated as a free left ${\cal F}({\mathbb C}_q^{2\vert1})$-module by the differentials
${\sf d}x_1$, ${\sf d}x_2$, ${\sf d}\theta$. We assume that the possible commutation relations between the generators and their first order differentials are of the form
\begin{align} \label{4.2}
{\sf d}v_i\cdot u_j = \sum_{k,l=1}^3 C_{ij}^{kl} \, u_k\cdot {\sf d}v_l
\end{align}
where the constants $C_{ij}^{kl}$ are possibly depending on $q$. That is, let the right ${\cal F}({\mathbb C}_q^{2\vert1})$-module structure on $\Gamma^1$ be completely defined by (\ref{4.2}). The following two lemmas will help the proof of Theorem \ref{thm4.1} below.

The existence of $\Delta_R$ depends on the commutation relations between $x_1:=x$, $x_2:=y$, $\theta$ and their differentials.

\begin{lemma} \label{lem4.1} 
Application of $\Delta_R$ to the relations $(\ref{4.2})$ from the left and the right leaves ten constants in $(\ref{4.2})$.
\end{lemma}

\begin{proof}
Using the identity in (\ref{4.1}), we can write, for $i=1,2$
\begin{equation} \label{4.3}
\Delta_R({\sf d}x) = {\sf d}x\otimes x, \quad \Delta_R({\sf d}y) = {\sf d}y\otimes y, \quad \Delta_R({\sf d}\theta) = {\sf d}\theta\otimes x^{-1}y.
\end{equation}
Application of the operator $\Delta_R$ from the left and the right sides to the relations (\ref{4.2}) imposes the following restrictions on some constants as follows
\begin{align*}
C_{ii}^{kl} &=0, \quad kl\ne ii; \quad
C_{13}^{31} = -q, \quad C_{13}^{kl}=0; \quad C_{23}^{32} = -q, \quad C_{23}^{kl}=0.
\end{align*}
As a result, it remains ten constants in the relations (\ref{4.2}).
\end{proof}

\begin{lemma} \label{lem4.2} 
Consistency with the relations of the algebra ${\cal F}({\mathbb C}_q^{2\vert1})$ leaves only two constants in $(\ref{4.2})$.
\end{lemma}

If we combine the above two lemmas, we obtain the following theorem, in which the two arbitrary constants remain in (\ref{4.2}).

\begin{theorem} \label{thm4.1} 
There exist two right-covariant ${\mathbb Z}_2$-graded first order differential calculi $\Gamma^1({\mathbb C}_q^{2\vert1})$ over the ${\mathbb Z}_2$-graded Hopf algebra
${\cal F}({\mathbb C}_q^{2\vert1})$ with respect to itself such that $\{{\sf d}x,{\sf d}y,{\sf d}\theta\}$ is a free right
${\cal F}({\mathbb C}_q^{2\vert1})$-module basis of  $\Gamma^1({\mathbb C}_q^{2\vert1})$. The bimodule structures for these calculi are determined by the relations
\begin{align} \label{4.4}
{\sf d}x\cdot x &= r x\cdot{\sf d}x, & {\sf d}\theta\cdot y &= q^{-1}s y\cdot{\sf d}\theta, \nonumber\\
{\sf d}y\cdot x &= r x\cdot {\sf d}y, & {\sf d}x\cdot\theta &= - q \theta\cdot{\sf d}x + (r-1)x\cdot{\sf d}\theta, \nonumber\\
{\sf d}\theta\cdot x &= q^{-1}r x\cdot{\sf d}\theta, & {\sf d}y\cdot\theta &= -q \theta\cdot{\sf d}y + (s-1) y\cdot{\sf d}\theta, \\
{\sf d}x\cdot y &= y\cdot{\sf d}x + (r-1) x\cdot{\sf d}y, & {\sf d}\theta\cdot\theta &= \theta\cdot{\sf d}\theta, \nonumber\\
{\sf d}y\cdot y &= s\,y\cdot{\sf d}y, \nonumber
\end{align}
for $C_{21}^{21}=0$, $C_{31}^{13}=rq^{-1}$, $C_{32}^{23}=sq^{-1}$, and
\begin{align} \label{4.5}
{\sf d}x\cdot x &= r x\cdot{\sf d}x, & {\sf d}\theta\cdot y &= q^{-1}s y\cdot{\sf d}\theta, \nonumber\\
{\sf d}y\cdot x &= x\cdot{\sf d}y+(s-1)y\cdot{\sf d}x, & {\sf d}x\cdot\theta &=-q\theta\cdot{\sf d}x+(r-1) x\cdot{\sf d}\theta, \nonumber\\
{\sf d}\theta\cdot x &= q^{-1}r\, x\cdot{\sf d}\theta, & {\sf d}y\cdot\theta &= -q\theta\cdot{\sf d}y+ (s-1) y\cdot{\sf d}\theta, \\
{\sf d}x\cdot y &= s\,y\cdot{\sf d}x, & {\sf d}\theta\cdot\theta &= \theta\cdot{\sf d}\theta, \nonumber\\
{\sf d}y\cdot y &= s\,y\cdot{\sf d}y, \nonumber
\end{align}
for $C_{12}^{12}=0$, $C_{31}^{13}=rq^{-1}$, $C_{32}^{23}=sq^{-1}$,
where $C_{11}^{11}=r$ and $C_{22}^{22}=s$.
\end{theorem}

\begin{remark} \label{rem4.1} 

We will see in Section 4.2 that relations $(\ref{4.5})$ are not consistent with higher order differential calculus (see, Remark \ref{rem4.6}).
Another choice of constants $C_{31}^{13}$ and $C_{32}^{23}$ for both $C_{21}^{21}=0$ and $C_{12}^{12}=0$ is $q$, but the resulting calculi are also inconsistent with higher order calculus.
\end{remark}

\begin{theorem} \label{thm4.2} 
There exists a map $\tau:{\cal F}({\mathbb C}_q^{2\vert1})\to M_3({\cal F}({\mathbb C}_q^{2\vert1}))$ that leads to relations $(\ref{4.4})$.
\end{theorem}

\begin{proof}
Let us consider right module structure of ${\cal F}({\mathbb C}_q^{2\vert1})$-bimodule $\Gamma^1({\mathbb C}_q^{2\vert1})$. We know that the right product
${\sf d}v \mapsto {\sf d}v\cdot u$ is an endomorphism of the left module $\Gamma^1({\mathbb C}_q^{2\vert1})$. Ring of all endomorphisms of any free module of rank 3 is isomorphic to the ring of all $3\times3$-matrices. Since $\{{\sf d}x, {\sf d}y, {\sf d}\theta\}$ is the homogeneous basis of $\Gamma^1({\mathbb C}_q^{2\vert1})$,  we can define a map $\tau$ by the formulas
\begin{equation} \label{4.6}
{\sf d}x_i\cdot u = \sum_{j=1}^3 (-1)^{p(u)p({\sf d}x_j)} \tau_{ij}(u)\cdot{\sf d}x_j
\end{equation}
for all $u\in {\cal F}({\mathbb C}_q^{2\vert1})$ and $x_1=x$, $x_2=y$ and $x_3=\theta$. Then, one can see that the relations (\ref{4.6}) equivalent to the relations (\ref{4.4}), where
\begin{align} \label{4.7}
\tau(x) &= \begin{pmatrix} r\,x & 0 & 0 \\ 0 & r\,x & 0 \\ 0 & 0 & q^{-1}r\,x \end{pmatrix}, \quad
\tau(y) = \begin{pmatrix} y & (r-1) x & 0 \\ 0 & s\,y & 0 \\ 0 & 0 & q^{-1}s\,y \end{pmatrix}, \nonumber\\
\tau(\theta) &= \begin{pmatrix} q\theta & 0 & (r-1) x \\ 0 & q\theta & (s-1) y \\ 0 & 0 & \theta \end{pmatrix}.
\end{align}
Here $r$ and $s$ are non-zero arbitrary complex numbers.
\end{proof}

\begin{theorem} \label{thm4.3} 
The map $\tau$ is a ${\mathbb Z}_2$-graded ${\mathbb C}$-linear homomorphism such that
\begin{equation} \label{4.8}
\tau_{ij}(uv) = \sum_{k=1}^3 (-1)^{p(u)[p(x_k)-p(x_j)]} \tau_{ik}(u) \tau_{kj}(v), \quad \forall u,v\in {\cal F}({\mathbb C}_q^{2\vert1})
\end{equation}
where $x_1=x$, $x_2=y$ and $x_3=\theta$.
\end{theorem}

\begin{proof}
The desired equality follows from the identity $({\sf d}x_i\cdot u)\cdot v={\sf d}x_i\cdot(uv)$. Indeed, for the left side of this equality we have
\begin{align*}
({\sf d}x_i\cdot u)\cdot v
&= \sum_{j=1}^3 (-1)^{p(u)p({\sf d}x_j)} \, \tau_{ij}(u)\cdot {\sf d}x_j \cdot v \\
&=\sum_{j=1}^3 (-1)^{p(u)p({\sf d}x_j)} \, \tau_{ij}(u)\cdot \sum_{k=1}^3 (-1)^{p(v)p({\sf d}x_k)} \, \tau_{jk}(v)\cdot {\sf d}x_k \\
&= \sum_{j=1}^3 (-1)^{p(v)p({\sf d}x_j)}\left[\sum_{k=1}^3 (-1)^{p(u)p({\sf d}x_k)} \, \tau_{ik}(u) \tau_{kj}(v)\right]\cdot {\sf d}x_j
\end{align*}
In the other hand, using (\ref{4.6}), we can write
\begin{align*}
{\sf d}x_i\cdot(uv) &= \sum_{j=1}^3 (-1)^{p(uv)p({\sf d}x_j)} \, \tau_{ij}(uv) \cdot {\sf d}x_j.
\end{align*}
As a result, we arrive at the desired equality with using the fact $p({\sf d}x_j)=p(x_j)+1$ (mod 2).
\end{proof}

\begin{remark} \label{rem4.2} 
It is easy to see that the relations $(\ref{3.1})$ are preserved under the action of the operator $\tau$.
\end{remark}

\begin{remark} \label{rem4.3} 
In particular, the linear functionals $f:{\cal F}({\mathbb C}_q^{2\vert1})\longrightarrow{\rm M}_3({\mathbb C})$ defined by
\begin{equation*}
f(x) = \begin{pmatrix} r & 0 & 0 \\ 0 & r & 0 \\ 0 & 0 & q^{-1}r \end{pmatrix}, \,
f(y) = \begin{pmatrix} 1 & r-1 & 0 \\ 0 & s & 0 \\ 0 & 0 & q^{-1}s \end{pmatrix}, \,
f(\theta) = \begin{pmatrix} 0 & 0 & r-1 \\ 0 & 0 & s-1 \\ 0 & 0 & 0 \end{pmatrix}
\end{equation*}
which are derived from $\tau$, are a representation of ${\cal F}({\mathbb C}_q^{2\vert1})$ with entries in ${\mathbb C}$, that is, they preserve the relations $(\ref{3.1})$ in the sense $0=f(x) f(y) - f(y) f(x) = f(xy-yx)$, {\it etc}.
\end{remark}

\begin{corollary} \label{cor4.1} 
We can also define a map  $\sigma:{\cal F}({\mathbb C}_q^{2\vert1})\to M_3({\cal F}({\mathbb C}_q^{2\vert1}))$ by the formulas
\begin{equation} \label{4.9}
u\cdot{\sf d}x_j = (-1)^{p(u)p({\sf d}x_j)} \sum_{i=1}^3 {\sf d}x_i\cdot\sigma_{ij}(u)
\end{equation}
for all $u\in {\cal F}({\mathbb C}_q^{2\vert1})$ and $x_1=x$, $x_2=y$ and $x_3=\theta$. We have
\begin{align}\label{4.10}
\sigma(x) &= \begin{pmatrix} r^{-1}x & 0 & 0 \\ 0 & r^{-1}x & 0 \\ 0 & 0 & r^{-1}q\,x \end{pmatrix}, \quad
\sigma(y) = \begin{pmatrix} y & 0 & 0 \\ (r^{-1}-1)x & s^{-1}y & 0 \\ 0 & 0 & s^{-1}q y \end{pmatrix}, \nonumber\\
\sigma(\theta) &= \begin{pmatrix} q^{-1} \theta & 0 & 0 \\ 0 & q^{-1}\theta & 0\\ (r^{-1}-1)x & (s^{-1}-1)y & \theta \end{pmatrix}
\end{align}
for the relations $(\ref{4.4})$.
\end{corollary}

\begin{corollary} \label{cor4.2} 
The map $\sigma$ is a ${\mathbb Z}_2$-graded ${\mathbb C}$-linear homomorphism in the sense $(\ref{4.8})$.
\end{corollary}

\begin{remark} \label{rem4.4} 
It is easy to see from Corollary \ref{cor4.2} that the relations $(\ref{3.1})$ are preserved under the action of the operator $\sigma$.
\end{remark}

\begin{remark} \label{rem4.5} 
In particular, the linear functionals $g:{\cal F}({\mathbb C}_q^{2\vert1})\longrightarrow{\rm M}_3({\mathbb C})$ defined by
\begin{align*}
g(x) &= \begin{pmatrix} 1 & 0 & 0 \\ 0 & 1 & 0 \\ 0 & 0 & q \end{pmatrix}, \quad
g(y) =\begin{pmatrix} s & 0 & 0 \\ s(r^{-1}-1) & 1 & 0 \\ 0 & 0 & q \end{pmatrix}, \nonumber\\
g(\theta) &= \begin{pmatrix} 0 & 0 & 0 \\ 0 & 0 & 0\\ r^{-1}-1 & s^{-1}-1 & 0 \end{pmatrix}
\end{align*}
which are derived from $\sigma$, are a representation of ${\cal F}({\mathbb C}_q^{2\vert1})$ with entries in ${\mathbb C}$, that is, they preserve the relations $(\ref{3.1})$.
\end{remark}

\subsection{Higher order differential calculus on ${\cal F}({\mathbb C}_q^{2\vert1})$}
\label{Subsec:4.2}

To form the de Rham complex of a $C^\infty$-manifold, we extend the exterior differential to the algebra of differential forms. We will consider a differential calculus as the de Rham complex for the algebra ${\cal F}({\mathbb C}_q^{2\vert1})$.

Let us begin with the definition of the ${\mathbb Z}_2$-graded differential calculus. Let ${\cal A}$ be an arbitrary associative (in general, noncommutative) algebra.

\begin{definition} \label{def4.4} 
A ${\mathbb Z}_2$-graded differential calculus over ${\cal A}$ is a ${\mathbb Z}_2$-graded algebra $\Omega^\wedge=\bigoplus_{n=0}^\infty\Omega^{\wedge n}$ where $\Omega^{\wedge0}={\cal A}$ and the space $\Omega^{\wedge n}$ of $n$-forms are generated as ${\cal A}$-bimodules via the action of a
${\mathbb C}$-linear mapping ${\sf d}:\Omega^\wedge \longrightarrow\Omega^\wedge$ of degree 1 such that

$1.$\, ${\sf d}^2 = 0$,

$2.$\, ${\sf d}(\alpha\wedge\beta) = ({\sf d}\alpha)\wedge\beta + (-1)^{p(\alpha)} \, \alpha\wedge({\sf d}\beta)$ for $\alpha,\beta\in\Omega^{\wedge}$.

$3.$\, $\Omega^{\wedge n}={\rm Lin}\{a_0\cdot{\sf d}a_1\wedge\cdots\wedge {\sf d}a_n: \, a_0,a_1,\ldots,a_n\in{\cal A}\}$ for $n\in{\mathbb N}$.
\end{definition}

Differential forms of higher order are obtained by applying {\sf d} to 1-forms (and then also higher forms) using ${\sf d}^2 = 0$ and the graded Leibniz rule.

We now introduce the first order differentials of the generators of ${\cal F}({\mathbb C}_q^{2\vert1})$ as ${\sf d}x=\varphi_1$,
${\sf d}y=\varphi_2$ and ${\sf d}\theta=z$. Then the differential ${\sf d}$ is uniquely defined by the conditions in Definition \ref{def4.4}, and the commutation relations between the differentials have the form (see, Definition \ref{def3.2})
\begin{equation} \label{4.11}
\varphi_i\wedge\varphi_j = - q^{2(i-j)} \varphi_j\wedge\varphi_i, \quad z\wedge\varphi_i = q\, \varphi_i\wedge z
\end{equation}
for $i,j=1,2$.

\begin{remark} \label{rem4.6} 
(i) If we use the conditions 1 and 2 in Definition \ref{def4.4}, we see that the relations $(\ref{4.5})$ are not consistent with relations $(\ref{4.11})$. However, application of the operator {\sf d} from the left to the relations $(\ref{4.4})$ says that the arbitrary constants $r$ and $s$ for the first-order calculus must equal to $q^2$. Therefore, we will continue with the relations $(\ref{4.4})$.

(ii) Using the commutation relations (\ref{4.4}), it is possible to find an $R$-matrix  that obeys the ${\mathbb Z}_2$-graded Yang-Baxter equation. If it is assumed that an $R$-matrix is associated with the ${\mathbb Z}_2$-graded space ${\mathbb C}_q^{2\vert1}$, then we see that the relations (\ref{4.4})
can be expressed as
\begin{equation*}
x_j\cdot{\sf d}x_k = q^2 \sum_{m,n=1}^3 (-1)^{p(x_m)} \,\hat{R}_{jk}^{mn} \, x_m\cdot {\sf d}x_n.
\end{equation*}
Here the entries of the matrix $\hat{R}$ are
$\hat{R}^{11}_{11}=\hat{R}^{21}_{12}=\hat{R}^{22}_{22}=1$, $\hat{R}^{12}_{12}=\hat{R}^{13}_{13}=\hat{R}^{23}_{23}=1-q^{-2}$,
$\hat{R}^{12}_{21}=-\hat{R}_{33}^{33}=q^{-2}$, $\hat{R}^{13}_{31}=\hat{R}^{23}_{32}=\hat{R}^{31}_{13}=\hat{R}^{32}_{23}=q^{-1}$ , except the zero entries. The matrix $R$ is given by $R=P\hat{R}$ where $P$ is the ${\mathbb Z}_2$-graded permutation matrix. The matrix $R$ satisfies the graded Yang-Baxter equation $R_{12}R_{13}R_{23}=R_{23}R_{13}R_{12}$ where $R_{12}=R\otimes I_3$, $R_{23} = I_3\otimes R$ and $R_{13} = (P\otimes I_3) R_{23} (P\otimes I_3)$ with the 3x3 identity matrix $I_3$. The matrix $\hat{R}$ obeys braid relation $\hat{R}_{12}\hat{R}_{23}\hat{R}_{12} = \hat{R}_{23}\hat{R}_{12}\hat{R}_{23}$. We not give more details about the matrix $R$, here. We will do a detailed review in the next article.
\end{remark}

We will denote the ${\mathbb Z}_2$-graded differential algebra generated by the elements of the set $\{x,y,\theta,{\sf d}x,{\sf d}y,{\sf d}\theta\}$ and the relations $(\ref{4.4})$ and $(\ref{4.11})$ by $\Omega^\wedge$

\begin{theorem} \label{thm4.4} 
The ${\mathbb Z}_2$-graded calculus $(\Omega^\wedge,{\sf d})$ is a ${\mathbb Z}_2$-graded $*$-calculus for the ${\mathbb Z}_2$-graded $*$-algebra
${\cal O}({\mathbb R}_q^{2\vert1})$.
\end{theorem}

\begin{proof}
We must show that the relations (\ref{4.6}) are invariant under the involution. As it can be seen easily, the matrices $\tau$ are related to the matrices $\sigma$ by
\begin{equation*}
\tau_{ji}(u) = (-1)^{p(\sigma_{ij}(u))} [\sigma_{ij}(u)]^* , \qquad \forall u\in{\cal F}({\mathbb C}_q^{2\vert1})
\end{equation*}
where $Q\bar{Q}=1$ for all $Q\in\{q,r,s\}$. Now, since the involution of the left side of the equality in (\ref{4.9}) is
\begin{align*}
(u\cdot{\sf d}x_j)^* = (-1)^{p(u)p({\sf d}x_j)}\, ({\sf d}x_j)^*\cdot u^*  = -(-1)^{p(u)+[1+p(u)]p({\sf d}x_j)}\, {\sf d}x_j\cdot u,
\end{align*}
for all $u,x_j\in{\cal F}({\mathbb C}_q^{2\vert1})$, the involution of (\ref{4.9}) gives
\begin{align*}
(-1)^{p(u {\sf d}x_j)} {\sf d}x_j\cdot u = \sum_{i=1}^3 (-1)^{p(\tau_{ji}(u){\sf d}x_i)+p(\tau_{ji}(u))p({\sf d}x_i)} \tau_{ji}(u)\cdot{\sf d}x_i
\end{align*}
which is equivalent to the relations (\ref{4.6}).
\end{proof}

\subsection{The ${\mathbb Z}_2$-graded partial derivatives}
\label{Subsec:4.3}

In order to get commutation relations of the generators of ${\cal F}({\mathbb C}_q^{2\vert1})$ with derivatives, we first  introduce the partial derivatives of the generators of the algebra. Since $(\Omega^\wedge,{\sf d})$ is a right covariant differential calculus, there are uniquely determined elements
$\partial_k(u)\in {\cal F}({\mathbb C}_q^{2\vert1})$ such that
\begin{equation}\label{4.12}
{\sf d}u = {\sf d}x \,\partial_x(u) + {\sf d}y \,\partial_y(u) + {\sf d}\theta \,\partial_\theta(u)
\end{equation}
for any element $u$ in ${\cal F}({\mathbb C}_q^{2\vert1})$, where $\partial_x=\partial/\partial x$, {\it etc}. For consistency, the degree of the derivative $\partial_\theta$ should be 1.

\begin{definition} \label{def4.5} 
The linear mappings $\partial_x, \partial_y, \partial_\theta:{\cal F}({\mathbb C}_q^{2\vert1})\longrightarrow {\cal F}({\mathbb C}_q^{2\vert1})$ defined by $(\ref{4.12})$ are called the {\it partial derivatives} of the calculus $(\Omega^\wedge,{\sf d})$.
\end{definition}
We know that the partial derivatives of the calculus $(\Omega^\wedge,{\sf d})$ satisfy the property $\partial_i(x_j) = \delta_{ij}$.

The next theorem gives the relations of the elements of ${\cal F}({\mathbb C}_q^{2\vert1})$ with their partial derivatives.

\begin{theorem} \label{thm4.5} 
The relations between the generators of ${\cal F}({\mathbb C}_q^{2\vert1})$ and partial derivatives are as follows
\begin{align}\label{4.13}
\partial_x x & = 1 + q^{-2}x \partial_x, & \partial_y x &= q^{-2} x \partial_y, & \partial_\theta x & = q^{-1} x\partial_\theta, \nonumber\\
\partial_y y &= 1 + q^{-2}y \partial_y + (q^{-2}-1)x\partial_x, & \partial_x y &= y \partial_x,
  & \partial_\theta y &= q^{-1} y\partial_\theta, \\
\partial_\theta \theta &= 1 - \theta \partial_\theta + (q^{-2}-1)(x\partial_x+y\partial_y), & \partial_x \theta &= q^{-1}\theta \partial_x, & \partial_y \theta &= q^{-1}\theta \partial_y.\nonumber
\end{align}
\end{theorem}

\begin{proof}
It can be shown that partial derivatives $\partial_1:=\partial_x$, $\partial_2:=\partial_y$,  $\partial_3:=\partial_\theta$ and the homomorphism $\sigma$ are related by
\begin{equation}\label{4.14}
\partial_j(u\cdot v) = \partial_j(u)\cdot v + \sum_{i=1}^3 (-1)^{p(u) p(\partial_j)} \sigma_{ji}(u) \cdot \partial_i(v), \quad \forall u,v \in {\cal F}({\mathbb C}_q^{2\vert1}).
\end{equation}
Now inserting the matrices in (\ref{4.10}) into (\ref{4.14}) with the fact that $\partial_i(x_j) = \delta_{ij}$ will give the relations (\ref{4.13}).
\end{proof}

The proof of the following theorem follows from the fact that ${\sf d}^2f = 0$ for a differentiable function $f$.

\begin{theorem} \label{thm4.6} 
The partial derivatives satisfy the following commutation relations
\begin{equation} \label{4.15}
\partial_x \partial_y = q^{-2}\partial_y \partial_x, \quad \partial_x \partial_\theta = q^{-1}\partial_\theta \partial_x, \quad
 \partial_y \partial_\theta = q^{-1} \partial_\theta \partial_y, \quad \partial_\theta^2 = 0.
\end{equation}
\end{theorem}

\subsection{${\mathbb Z}_2$-graded right-invariant Maurer-Cartan 1-forms}
\label{Subsec:4.4}

In order to determine the most general commutation relations of elements of $\Omega^\wedge$ with elements of ${\cal F}({\mathbb C}_q^{2\vert1})$ we use a convenient basis of $\Omega^\wedge$. It consists of the quantum analogues of the Maurer-Cartan 1-forms defined by
\begin{equation} \label{4.16}
w_i:=w(u_i) = m({\sf d}\otimes S)\Delta(u_i), \quad u_i\in {\cal F}({\mathbb C}_q^{2\vert1}).
\end{equation}

\begin{lemma} \label{lem4.3} 
$(1)$ The 1-forms $w_i$ are right-invariant, i.e.
\begin{equation*}
\Delta_R(w_i) = w_i\otimes{\bf 1}.
\end{equation*}
$(2)$ The set ${\cal B}= \{w_i\}_{i=1,2,3}$ is a basis of $\Omega_{\rm inv}^\wedge$ as a ${\mathbb C}$-vector space.
\end{lemma}

The following theorem determines the relations between right-coinvariant forms and the elements of ${\cal F}({\mathbb C}_q^{2\vert1})$.

\begin{theorem} \label{thm4.7} 
Let $(\Omega^\wedge,\Delta_R)$ be a right-covariant bimodule over ${\cal F}({\mathbb C}_q^{2\vert1})$ and $\{w_k\}$ be a basis in the vector space of all right-invariant elements of $\Omega^\wedge$. Then there exist a linear map
$\mu:{\cal F}({\mathbb C}_q^{2\vert1})\longrightarrow M_3({\cal F}({\mathbb C}_q^{2\vert1}))$ such that
\begin{align} \label{4.17}
w_i \cdot u &= \sum_{j=1}^3 \, (-1)^{p(u)p(w_i)} \, \mu_{ij}(u)\cdot w_j
\end{align}
for all $u\in {\cal F}({\mathbb C}_q^{2\vert1})$.
\end{theorem}

\begin{proof}
Using (\ref{4.16}) and (\ref{4.4}), we find the matrices $\mu(u)$ for $u\in {\cal F}({\mathbb C}_q^{2\vert1})$ as
\begin{align} \label{4.18}
\mu(x) &= \begin{pmatrix} q^2x & 0 & 0 \\ 0 & q^2x & 0\\ 0 & 0 & qx \end{pmatrix}, \quad
  \mu(y) = \begin{pmatrix} y & (q^2-1)y & 0\\ 0 & q^2y & 0\\ 0 & 0 & qy \end{pmatrix}, \nonumber\\
\mu(\theta) &= \begin{pmatrix} \theta & 0 & (q^{-2}-1)x^{-1}y \\ 0 & \theta & (q^{-2}-1)x^{-1}y \\ 0 & 0 & \theta \end{pmatrix}.
\end{align}
Inserting these matrices in (\ref{4.17}) leads to relations of Maurer-Cartan 1-forms with the generators of the algebra ${\cal F}({\mathbb C}_q^{2\vert1})$.
\end{proof}

\begin{corollary} \label{cor4.4} 
The map $\mu$ is ${\mathbb C}$-linear ${\mathbb Z}_2$-graded homomorphism such that
\begin{equation*}
\mu_{ij}(u\cdot v) = \sum_{k=1}^4 (-1)^{p(v)[p(x_k)-p(x_i)]} \,\mu_{ik}(u)\mu_{kj}(v), \qquad (i,j=1,2,3)
\end{equation*}
where $x_i,x_k,u,v\in {\cal F}({\mathbb C}_q^{2\vert1})$.
\end{corollary}

\begin{proof}
Using (\ref{4.17}), we can write
\begin{align*}
w_i\cdot(uv) &= \sum_{j=1}^3 (-1)^{p(uv)p(w_i)} \, \mu_{ij}(uv) w_j.
\end{align*}
In the other hand, we have
\begin{align*}
(w_i\cdot u)v
&= \left[\sum_{j=1}^3 (-1)^{p(u)p(w_i)} \, \mu_{ij}(u)\cdot w_j\right]\cdot v \\
&=\sum_{j=1}^3 (-1)^{p(u)p(w_i)} \, \mu_{ij}(u)\cdot\left[\sum_{k=1}^3 (-1)^{p(v)p(w_j)} \, \mu_{jk}(v)\cdot w_k\right] \\
&= \sum_{j=1}^3 (-1)^{p(u)p(w_i)}\left[\sum_{k=1}^3 (-1)^{p(v)p(w_k)} \, \mu_{ik}(u) \mu_{kj}(v)\right]\cdot w_j.
\end{align*}
As a result, we arrive at the desired equality with the identity $w_i\cdot(uv) = (w_i\cdot u)v$. Here we used the fact that $p(w_i)=p(x_i)+1$ (mod 2).
\end{proof}

We know that from \cite{21}, if $(\Omega,\Delta_R)$ is a right-covariant bimodule over a Hopf algebra ${\cal H}$ and $\{w_k\}$ is a basis in the vector space of all right-invariant elements of $\Omega$ then there exist linear functionals $F_{ij}$ such that
\begin{equation} \label{4.19}
w_i \cdot u = \sum_j \,m\circ(F_{ij}\otimes {\rm id})\Delta(u)\, w_j
\end{equation}
for all $u\in {\cal H}$. The ${\mathbb Z}_2$-graded form of this formula is expressed in the following corollary.

\begin{corollary} \label{cor4.5} 
Let $(\Omega^\wedge,\Delta_R)$ be a right-covariant bimodule over ${\cal F}({\mathbb C}_q^{2\vert1})$ and $\{w_k\}$ be a basis in the vector space of all right-invariant elements of $\Omega^\wedge$. Then there exist linear map $F: {\cal F}({\mathbb C}_q^{2\vert1})\longrightarrow {\rm M}_3({\mathbb C})$ such that
\begin{equation} \label{4.20}
w_i \cdot u = \sum_{j=1}^3 (-1)^{p(u)p(w_i)} \, m\circ(F_{ij}\otimes {\rm id})\Delta(u)\, w_j
\end{equation}
for all $u\in {\cal F}({\mathbb C}_q^{2\vert1})$.
\end{corollary}

\begin{proof}
Indeed, substituting the matrices
\begin{equation} \label{4.21}
F(x) = \left(\begin{matrix} q^2 & 0 & 0 \\ 0 & q^2 & 0\\ 0 & 0 & q \end{matrix}\right), \,
F(y) = \left(\begin{matrix} 1 & q^2-1 & 0\\ 0 & q^2 & 0\\ 0 & 0 & q \end{matrix}\right), \,
F(\theta) =\left(\begin{matrix} 0 & 0 & q^{-2}-1 \\ 0 & 0 & q^{-2}-1 \\ 0 & 0 & 0 \end{matrix}\right)
\end{equation}
in (\ref{4.20}) we arrive at the relations (\ref{4.17}).
\end{proof}

\begin{remark} \label{rem4.7} 
It is easy to see that the matrices $f(x_i)$ in (\ref{4.21}) are a representation of ${\cal F}({\mathbb C}_q^{2\vert1})$ with entries in ${\mathbb C}$, that is, they satisfy the relations $(\ref{3.1})$ in the sense $(\ref{4.8})$.
\end{remark}

\begin{corollary} \label{cor4.6} 
There is a kind of inverse formula of $(\ref{4.17})$ as follows
\begin{align*}
u \cdot w_i &= \sum_{j=1}^3 (-1)^{p(u)p(w_j)} \, w_j \cdot\tilde{\mu}_{ji}(u), \quad \forall u\in {\cal F}({\mathbb C}_q^{2\vert1})
\end{align*}
where $\tilde{\mu}: {\cal F}({\mathbb C}_q^{2\vert1})\longrightarrow {\rm M}_3({\mathbb C}_q^{2\vert1}))$ is a ${\mathbb Z}_2$-graded algebra homomorphism defined by
\begin{align} \label{4.22}
\tilde{\mu}(x) &= \left(\begin{matrix} q^{-2}x & 0 & 0 \\ 0 & q^{-2}x & 0\\ 0 & 0 & q^{-1}x \end{matrix}\right), \quad
  \tilde{\mu}(y) = \left(\begin{matrix} y & 0 & 0\\ (q^{-2}-1)y & q^{-2}y & 0\\ 0 & 0 & q^{-1}y \end{matrix}\right), \nonumber\\
\tilde{\mu}(\theta) &=\left(\begin{matrix} \theta & 0 & 0 \\ 0 & \theta & 0 \\ (1-q^{-2})x^{-1}y & (1-q^{-2})x^{-1}y & \theta \end{matrix}\right).
\end{align}
\end{corollary}

\begin{remark} \label{rem4.9} 
The operators  $\mu$ and $\tilde{\mu}$ are related by the formulae
\begin{equation*}
\sum_{k=1}^3 \tilde{\mu}_{jk}(\mu_{ik}(u)) = u \, \delta_{ij}, \quad \forall u\in {\cal F}({\mathbb C}_q^{2\vert1}).
\end{equation*}
\end{remark}

Using (\ref{4.4}) and (\ref{4.17}) we find the commutation rules of the forms as follows.

\begin{theorem} \label{thm4.8} 
The commutation rules of the forms are as follows
\begin{align} \label{4.23}
w_1\wedge w_1 &= 0, & w_1\wedge w_2 &= - w_2\wedge w_1, \quad w_2 \wedge w_3 = w_3\wedge w_2, \nonumber\\
w_2\wedge w_2 &= 0, & w_1\wedge w_3 &= q^{-2} \, w_3\wedge w_1 + (1-q^{-2}) \, w_3\wedge w_2.
\end{align}
\end{theorem}

\begin{theorem} \label{thmm4.9} 
The Maurer-Cartan equations for $\Omega^\wedge$ have the form
\begin{equation} \label{4.24}
{\sf d}w_1 = 0, \quad {\sf d}w_2 = 0, \quad {\sf d}w_3 = q^2w_3\wedge(w_1-w_2).
\end{equation}
\end{theorem}

\subsection{A special automorphism group}
\label{Subsec:4.5}

Here we will introduce an upper triangular matrix with six entries, and without going into much detail we will emphasize that such matrices form a ${\mathbb Z}_2$-graded quantum group (or supergroup) as a ${\mathbb Z}_2$-graded Hopf algebra and that the algebras ${\cal O}({\mathbb C}_q^{2\vert1})$ and ${\cal O}({\mathbb C}_q^{1\vert2}):=\Lambda({\mathbb C}_q^{2\vert1})$ are the left co-module algebras of this quantum supergroup. Moreover, we will see that the calculus given in Section 4 is the left-covariant under this group.

Let ${\cal O}({\rm M}^u(2\vert1))$ be defined as the polynomial algebra ${\mathbb C}[a,b,c,d,\alpha,\beta]$ where $a,b,c$, $d,\alpha,\beta$ are elements of an associative ${\mathbb Z}_2$-graded algebra ${\cal A}$. We assume the generators $a$, $b$, $c$, $d$ are of degree 0 and the generators $\alpha,\beta$ are of degree 1, and write a point $(a,b,c,d,\alpha,\beta)$ of ${\cal O}({\rm M}^u(2\vert1))$ in the matrix form, as a supermatrix,
\begin{equation*}
T = \begin{pmatrix} a & bd & \alpha c \\ 0 & b & \beta c\\ 0 & 0 & c\end{pmatrix} =(t_{ij}).
\end{equation*}

If the generators of the algebras ${\cal O}({\mathbb C}_q^{2\vert1})$ and ${\cal O}({\mathbb C}_q^{1\vert2}))$ are written in vector form as ${\bf x}=(x_1,x_2,\theta)^t$ and $\hat{\bf x}=(\varphi_1,\varphi_2,z)^t$, and if the vectors ${\bf x}$ and $\hat{\bf x}$ are multiplied by the matrix $T$ from the left, the resulting elements are again expected to satisfy the relations (\ref{3.1}) and (\ref{3.2}). For this to happen, it is sufficient to assume that the components of ${\bf x}$ and $\hat{\bf x}$ are supercommutative with matrix elements of the matrix $T$. Consequently, we have the following theorem whose proof can be done by direct computations using the relations (\ref{3.1}) and (\ref{3.2}).

\begin{theorem} \label{thm4.9} 
The entries of the supervectors ${\bf x}'$ and $\hat{\bf x}'$ satisfy $(\ref{3.1})$ and $(\ref{3.2})$, respectively, if and only if the entries of $T$ fulfill the relations
\begin{align*}
ab &= ba, & ac &= ca, & bc &= cb, & \alpha^2 &= 0, \nonumber\\
ad &= q^2da, & bd &= db, & cd &= dc, & \beta^2 &= 0, \nonumber\\
a\alpha &= q\alpha a, & b\alpha &= q \alpha b, & c\alpha &= q \alpha c, & \beta\alpha &= -q^2 \alpha\beta, \\
a\beta &= q^{-1}\beta a, & b\beta &= q \beta b, & c\beta &= q \beta c, \nonumber\\
\alpha d &= d\alpha, & \beta d &= q^2 d\beta + (1-q^2) \alpha, \nonumber
\end{align*}
\end{theorem}
The algebra ${\cal O}({\rm M}_q^u(2\vert1))$ is the quotient of the free algebra ${\mathbb C}\langle a,b,c,d,\alpha,\beta\rangle$ by the two-sided ideal $J_q$ generated by the relations given in Theorem \ref{thm4.9}. The superdeterminant of the matrix $T$ is given by ${\rm sdet}(T)=abc^{-1}$, where the formal inverse of $c$ exists.
Let ${\cal O}({\rm GL}_q^u(2\vert1))$ be the quotient of the algebra ${\cal O}({\rm M}_q^u(2\vert1))$ by the two-sided ideal generated by the element
$t\cdot{\rm sdet}(T)-1$. The ${\mathbb Z}_2$-graded Hopf algebra structure of ${\cal O}({\rm GL}_q^u(2\vert1))$ is given as usual in the following theorem.

\begin{theorem} \label{thm4.10} 
(i) There exists a unique ${\mathbb Z}_2$-graded Hopf algebra structure on the algebra ${\cal O}({\rm GL}_q^u(2\vert1))$ with co-maps $\Delta$, $\epsilon$ and $S$ such that
\begin{eqnarray*}
\Delta(t_{ij}) = \sum_{k=1}^3 t_{ik} \otimes t_{kj}, \quad \epsilon(t_{ij}) = \delta_{ij}, \quad S(T)=T^{-1}.
\end{eqnarray*}

\noindent(ii) The ${\mathbb Z}_2$-graded algebras ${\cal O}({\mathbb C}_q^{2\vert1})$ and ${\cal O}({\mathbb C}_q^{1\vert2})$ are left ${\cal O}({\rm GL}_q^u(2\vert1))$-comodule algebras.

\noindent(iii) Differential calculus $\Omega^\wedge({\mathbb C}_q^{2\vert1})$ over the algebra ${\cal O}({\mathbb C}_q^{2\vert1})$ is left-covariant with respect to the ${\mathbb Z}_2$-graded Hopf algebra ${\cal O}({\rm GL}_q^u(2\vert1))$.
\end{theorem}

Therefore the quantum supergroup ${\rm GL}_q^u(2\vert1)$ can be considered as a quantum automorphism supergroup of a pair of quantum linear ${\mathbb Z}_2$-graded spaces
${\cal O}({\mathbb C}_q^{2\vert1})$ and ${\cal O}({\mathbb C}_q^{1\vert2})$ defined by the relations (\ref{3.1}) and (\ref{3.2}), respectively. We will do a detailed review of such quantum supergroups in the next paper.

\subsection{${\mathbb Z}_2$-graded Quantum Weyl algebra}
\label{Sec:4.6}

In this section, we will define the involution for partial derivatives and give some related topics.

In Section 3, we considered ${\cal O}({\mathbb C}_q^{2\vert1})$ in terms of generators and relations and define ${\cal O}({\mathbb C}_q^{2\vert1})$ as the quotient $
{\mathbb C}\langle x_1,x_2,\theta \rangle/I_q$. In the limit $q\to1$ the relations (\ref{3.1}) reduce to the defining relations of the ${\mathbb Z}_2$-graded algebra
${\cal O}({\mathbb C}^{2\vert1})$. We also noted that the set $\{x^my^n\theta^k: \, m,n\in{\mathbb N}_0, k=0,1\}$ forms a vector space basis of ${\cal O}({\mathbb C}_q^{2\vert1})$. Below we give the definition of a dual algebra.

\begin{definition} \label{def4.6} 
Let ${\mathbb C}\langle \partial_x,\partial_y,\partial_\theta \rangle$ be a free ${\mathbb Z}_2$-graded unital algebra generated by the elements $\partial_x$, $\partial_y$, $\partial_\theta$. Then we define the ${\mathbb Z}_2$-graded algebra ${\cal O}(\partial{\mathbb C}_q^{2\vert1})$ as
$${\cal O}(\partial{\mathbb C}_q^{2\vert1}) = {\mathbb C}\langle \partial_x,\partial_y,\partial_\theta \rangle/I_\partial$$
where $I_\partial$ is spanned by quadratic expressions of the form $\partial_x \partial_y-q^{-2}\partial_y \partial_x$, $\partial_x \partial_\theta-q^{-1}\partial_\theta \partial_x$, $\partial_y \partial_\theta-q^{-1} \partial_\theta \partial_y$, $\partial_\theta^2$.
\end{definition}
It can be shown that the set of monomials $\{\partial_x^m \partial_y^n \partial_\theta^k: \, m,n\in{\mathbb N}_0, k=0,1\}$ forms a vector space basis of
${\cal O}(\partial{\mathbb C}_q^{2\vert1})$. In the case $q=1$ the relations (\ref{4.13}) and (\ref{4.15}) together with (\ref{3.1}) reduce to the defining relations of the
${\mathbb Z}_2$-graded Weyl algebra denoted by ${\cal W}(2\vert1)$. This motivates following definition.

\begin{definition} \label{def4.7} 
The ${\mathbb Z}_2$-graded quantum Weyl algebra ${\cal W}_q(2\vert1)$ is the unital algebra generated by $x$, $y$, $\theta$ and $\partial_x$, $\partial_y$, $\partial_\theta$ and the relations $(\ref{3.1})$, $(\ref{4.13})$ and $(\ref{4.15})$.
\end{definition}

We mention a few algebraic properties of the algebra ${\cal W}_q(2\vert1)$.

{\bf 1.} We know that the relations (\ref{3.1}) generate $I_q$ and (\ref{4.15}) generate $I_\partial$. So, there are algebra maps
${\cal O}({\mathbb C}_q^{2\vert1})\longrightarrow {\cal W}_q(2\vert1)$ and ${\cal O}(\partial{\mathbb C}_q^{2\vert1})\longrightarrow {\cal W}_q(2\vert1)$.
Since the linear map of ${\cal O}({\mathbb C}_q^{2\vert1})\otimes {\cal O}(\partial{\mathbb C}_q^{2\vert1})$ to ${\cal W}_q(2\vert1)$ defined by $u\otimes v \mapsto u\cdot v$ is a vector space isomorphism, one can consider ${\cal O}({\mathbb C}_q^{2\vert1})$ and ${\cal O}(\partial{\mathbb C}_q^{2\vert1})$ as subalgebras of ${\cal W}_q(2\vert1)$ and the set $\{x^my^n\theta^k\partial_x^{m'} \partial_y^{n'} \partial_\theta^{k'}: \, m,m',n,n'\in{\mathbb N}_0, k,k'=0,1\}$ is a vector space basis of
${\cal W}_q(2\vert1)$.

{\bf 2.} By Theorem \ref{thm4.4}, the first order differential calculus $\Omega^1$ given in Section 4.2 is a $*$-calculus for the $*$-algebra
${\cal O}({\mathbb C}_q^{2\vert1})$. The involution of $\Omega^1$ induces the involution for the partial differential operators:

\begin{lemma} \label{lem4.4} 
The algebra ${\cal O}(\partial{\mathbb C}_q^{2\vert1})$ supplied with the involutions determined by
\begin{equation} \label{4.25}
\partial_x^* = -q^{-2}\partial_x, \quad \partial_y^* = -q^{-4} \partial_y, \quad \partial_\theta^* = q^{-4} \partial_\theta,
\end{equation}
becomes a $*$-algebra ${\cal O}(\partial{\mathbb R}_q^{2\vert1})$, where $|q|=1$.
\end{lemma}

\begin{theorem} \label{thm4.11} 
If $|q|=1$ then the algebra ${\cal W}_q(2\vert1)$ is a $*$-algebra.
\end{theorem}

\begin{proof}
Using the assumption $|q|=1$ we check that the relations (\ref{3.1}), (\ref{4.13}) and (\ref{4.15}) are preserved under the involution. We know from Theorem \ref{thm3.2} and Lemma \ref{lem4.4} that the relations (\ref{3.1}) and (\ref{4.15}) are preserved under the involution. It remains to show that the relations (\ref{4.13}) are preserved under the involution, which is not difficult.
\end{proof}

{\bf 3.} As an element of the ${\mathbb Z}_2$-graded Weyl algebra ${\cal W}_q(2\vert1)$, the element
\begin{equation} \label{4.26}
{\sf D}=x\partial_x+y\partial_y+\theta\partial_\theta
\end{equation}
is called the {\it Euler derivation}. The action of ${\sf D}$ on ${\cal O}({\mathbb C}_q^{2\vert1})$ involves the $q$-integers $[n]_{q^{-2}}$. We set $[n]_{q^{-2}}=(q^{-2n}-1)/(q^{-2}-1)$ where $q^2\ne1$.

\begin{theorem} \label{thm4.12} 
(i) If $f_n\in {\cal O}({\mathbb C}_q^{2\vert1})$ and $g_n\in {\cal O}(\partial{\mathbb C}_q^{2\vert1})$ are homogeneous elements of order $n$, then we have the identities
\begin{equation} \label{4.27}
{\sf D}f_n = [n]_{q^{-2}} \,f_n + q^{-2n} \,f_n {\sf D} \quad {\rm and} \quad
  g_n {\sf D} = [n]_{q^{-2}} \,g_n + q^{-2n} \,{\sf D}\, g_n.
\end{equation}

\noindent(ii) The operator ${\sf D}_{q^{-2}}$ is a $q^{-2}$-derivative.

\noindent(iii) The element ${\sf E} = {\bf 1} + (q^{-2} - 1){\sf D}$ is a non-zero normal element in ${\cal W}_q(2\vert1)$.
\end{theorem}

\begin{proof}
(i) We know that every polynomial of order $n$ is a linear combination of monomials of order $n$ and so we may assume $f_n$ is a monomial. So, we first assume that  $f=xy\theta$. Then, we have
\begin{align*}
{\sf D}f
&= (x\partial_x+y\partial_y+\theta\partial_\theta) xy\theta = \left[x + q^{-2} x(x\partial_x+y\partial_y+\theta\partial_\theta)\right] y\theta\\
&=  xy\theta + q^{-2}x\left[y + q^{-2} y(x\partial_x+y\partial_y+\theta\partial_\theta)\right] \theta\\
&=  (1 + q^{-2}) xy\theta + q^{-4}xy\theta (x\partial_x+y\partial_y+\theta\partial_\theta) = (1 + q^{-2}) f + q^{-4} f {\sf D}.
\end{align*}
Now we get the first identity in (\ref{4.27}) by induction over $n$ by taking $f_n=x^my^{m'}\theta$, where $n=m+m'+1$. The proof of the second identity in (\ref{4.27}) is similar.

(ii) Let $q\ne1$. In terms of the ${\cal W}_q(2\vert1)$-module structure for ${\cal O}({\mathbb C}_q^{2\vert1})$, the part (i) implies that if $f_n$ is a homogeneous element of order $n$ in ${\cal O}({\mathbb C}_q^{2\vert1})$ then
\begin{align*}
\partial_x(x^my^{m'}\theta) &= q^m y^{m'}\theta \partial_x(x^m) = q^{-2m'}\cdot\frac{q^{-2m}-1}{q^{-2}-1} \, x^{m-1}y^{m'}\theta, \\
\partial_y(x^my^{m'}\theta) &= \frac{q^{-2m'}-1}{q^{-2}-1} \, x^my^{m'-1}\theta, \quad \partial_\theta(x^my^{m'}\theta) = q^{-2(m+m')} x^my^{m'}
\end{align*}
where $n=m+m'+1$, so that ${\sf D}(f_n)=[n]_{q^{-2}} \,f_n$. In particular,
\begin{equation*}
{\sf D}(f_nf_k) = [n+k]_{q^{-2}} \,f_nf_k = {\sf D}(f_n) f_k + q^{-2n}f_n{\sf D}(f_k).
\end{equation*}

(iii) It is easy to check that ${\sf E}u=q^{-2}u{\sf E}$ and $\partial_u{\sf E}=q^{-2}{\sf E}\partial_u$ for all $u\in {\cal F}({\mathbb C}_q^{2\vert1})$. We have also
${\sf E}^*=q^{-2}{\sf E}$.
\end{proof}

\section{${\mathbb Z}_2$-graded Hopf-Lie algebras}
\label{Sec:5}

In this section, we will derive three ${\mathbb Z}_2$-graded Lie algebras from the generators of the ${\mathbb Z}_2$-graded algebra ${\cal F}({\mathbb C}_q^{2\vert1})$, the calculus, and duality. As a result we will see that all three ${\mathbb Z}_2$-graded Lie algebras are equivalent with appropriate choices.

\subsection{A ${\mathbb Z}_2$-graded Lie algebra generated from ${\cal F}({\mathbb C}_q^{2\vert1})$}
\label{Subsec:5.1}

The generators $x$ and $y$ of the Hopf algebra ${\cal F}({\mathbb C}_q^{2\vert1})$ are group-like whereas for $\theta$ we have defined
\begin{equation} \label{5.1}
\Delta(\theta) = \theta \otimes x^{-1}y + {\bf 1}\otimes\theta, \quad S(\theta) = - q^{-1} x\theta y^{-1}.
\end{equation}
However, we know that
\begin{equation} \label{5.2}
\Delta(\theta) = \theta \otimes {\bf 1} + {\bf 1} \otimes \theta, \quad S(\theta) = - \theta
\end{equation}
for an odd element $\theta\in {\cal L}$ of universal enveloping ${\mathbb Z}_2$-graded algebra ${\cal U}(g)$ for a ${\mathbb Z}_2$-graded Lie algebra ${\cal L}$. So the definition (\ref{5.1}) is a twisting with $x$ and $y$ of (\ref{5.2}). Furthermore, if we consider the following correspondence between the generators of ${\cal F}({\mathbb C}_q^{2\vert1})$ and the ${\mathbb Z}_2$-graded Lie algebra generators
\begin{equation*}
x \longleftrightarrow q^{\tilde{x}}, \quad y \longleftrightarrow q^{\tilde{y}}, \quad \theta \longleftrightarrow \Phi
\end{equation*}
then we obtain a ${\mathbb Z}_2$-graded quantum Lie algebra ${\cal L}$, which is best stated in terms of $M=(\tilde{x}+\tilde{y})/h$ and $N=\tilde{x}-\tilde{y}$, for $\ln q=h\ne0$
\begin{equation} \label{5.3}
[M,N] = 0, \quad [M,\Phi] = 2\Phi, \quad [N,\Phi] = 0, \quad \Phi^2=0,
\end{equation}
where $[,]$ denotes Lie brackets. In this basis, the ${\mathbb Z}_2$-graded Hopf algebra structure takes the following form
\begin{align} \label{5.4}
\Delta(M) &= M\otimes{\bf 1} + {\bf 1}\otimes M, & \epsilon(M) &= 0, & S(M) &= -M, \nonumber\\
\Delta(N) &= N\otimes{\bf 1} + {\bf 1}\otimes N, & \epsilon(N) &= 0, & S(N) &= -N, \\
\Delta(\Phi) &= \Phi\otimes q^{-N} + {\bf 1}\otimes\Phi, & \epsilon(\Phi) &= 0, & S(\Phi) &= -\Phi q^{N}.\nonumber
\end{align}

\subsection{${\mathbb Z}_2$-graded Quantum Lie algebra}
\label{Sec:5.2}

Usually a given ${\mathbb Z}_2$-graded algebra is associated to a ${\mathbb Z}_2$-graded Lie algebra. The commutation relations of Maurer-Cartan forms allow us to construct the Lie algebra of the vector fields. In order to obtain the ${\mathbb Z}_2$-graded quantum Lie algebra of the vector fields we first write the differential {\sf d} in the form
\begin{equation*}
{\sf d}u = (\omega_1 T_1 + \omega_2 T_2 + \omega_3 \Theta)u, \quad u\in {\cal F}({\mathbb C}_q^{2\vert1})
\end{equation*}
where $T_1$, $T_2$ and $\Theta$ are the ${\mathbb Z}_2$-graded quantum Lie algebra generators (vector fields).

\begin{theorem} \label{thm5.1} 
The commutation relations of the ${\mathbb Z}_2$-graded quantum Lie algebra generators are in the form
\begin{equation*}
[T_1,T_2] = 0, \quad [T_1,\Theta]_{q^{-2}} = - q^2\Theta, \quad [T_2,\Theta] = q^2\Theta + (1-q^{-2})\Theta T_1, \quad \Theta^2=0
\end{equation*}
where $[A,B]_r = AB- (-1)^{P(A)p(B)} r BA$.
\end{theorem}

\begin{proof}
Considering an arbitrary function $f$ of the generators of ${\cal F}({\mathbb C}_q^{2\vert1})$ and using the fact that ${\sf d}^2 = 0$ with the relations (\ref{4.23}) and (\ref{4.24}) we obtain the desired relations.
\end{proof}

To better stated the relations given in the above theorem let us choose
\begin{equation*}
q^{-2H} = {\bf 1} + q^{-2}(1-q^{-2})T_1, \quad  q^{-2(H+K)} = {\bf 1} + (q^{-2}-1)(T_1 + T_2).
\end{equation*}
In this case, we have the relations
\begin{equation} \label{5.5}
[H,K] = 0, \quad [H,\Theta] = \Theta, \quad [K,\Theta] = -\Theta, \quad \Theta^2=0.
\end{equation}

We denote the ${\mathbb Z}_2$-graded Lie algebra generated by the unit {\bf 1} and the generators $H$, $K$, $\Theta$ modulo the relations (\ref{5.5}) by ${\cal L}$. We know that the generators $H$, $K$ and $\Theta$ are endowed with a natural coproduct and then ${\mathbb Z}_2$-graded Hopf algebra. To determine the action of co-maps on these vector fields, we first need a lemma that can be proved by the mathematical induction technique.

\begin{lemma} \label{lem5.1}  
The action of the 1-forms given in $(\ref{4.16})$ on the monomial $f=x^my^n\theta$ as follows
\begin{align} \label{5.6}
f w_1 &= -q^{-2m} w_1 f + q^{-2m}(1-q^{-2n}) w_2 f + (1-q^{-2})q^{-2(m+n)} w_3 x^{m-1}y^{n+1}, \nonumber\\
f w_2 &= - q^{-2(m+n)} w_2 f + (1-q^{-2})q^{-(m+n)} w_3 x^{k-1}y^{n+1}, \\
f w_3 &= q^{-(m+n)} w_3 f. \nonumber
\end{align}
\end{lemma}

\begin{theorem} \label{thm5.2} 
The ${\mathbb Z}_2$-graded Lie algebra ${\cal L}$ is a ${\mathbb Z}_2$-graded Hopf algebra. The coproduct, counit and coinverse on the ${\mathbb Z}_2$-graded Lie algebra ${\cal L}$ are as follows:
\begin{align} \label{5.7}
\Delta(H) &= H \otimes {\bf 1} + {\bf 1}\otimes H, & \epsilon(H) &= 0, & S(H) &= -H, \nonumber\\
\Delta(K) &= K \otimes {\bf 1} + {\bf 1}\otimes K, & \epsilon(K) &= 0, & S(K) &= -K, \\
\Delta(\Theta) &= \Theta\otimes{\bf 1} + q^{-(H+K)}\otimes\Theta & \epsilon(\Theta) &= 0, & S(\Theta) &= - q^{H+K}\Theta. \nonumber
\end{align}
\end{theorem}

\begin{proof}
The Leibniz rule may be expressed as follows
\begin{equation} \label{5.8}
{\sf d}(uv) = ((\omega_1 T_1 + \omega_2 T_2 + \omega_3 \Theta)u)v + u(\omega_1 T_1 + \omega_2 T_2 + \omega_3 \Theta)v.
\end{equation}
Inserting (\ref{5.6}) in (\ref{5.8}) and equating coefficients of the Maurer-Cartan 1-forms after neglecting some terms, we obtain the action of $\Delta$ on generators. The action of counit and coinvers on generators follows from Hopf algebra axioms in Definition \ref{def2.2}.
\end{proof}

\subsection{Dual pairings of ${\mathbb Z}_2$-graded Hopf algebras}
\label{Subsec:5.3}

Let us denote the ${\mathbb Z}_2$-graded Hopf algebra ${\cal F}({\mathbb C}_q^{2\vert1})$ by ${\cal A}$. Then its dual ${\cal U}\doteq {\cal A}'$ is a ${\mathbb Z}_2$-graded Hopf algebra as well. Using the coproduct $\Delta$ in ${\cal A}$, we can define a product in ${\cal U}$ and using the product in the ${\mathbb Z}_2$-graded Hopf algebra ${\cal A}$, we can define a coproduct in ${\cal U}$.

In this section, to get the dual of the ${\mathbb Z}_2$-graded Hopf algebra ${\cal A}$ defined in Section 3, we have applied to ${\cal A}={\cal F}({\mathbb C}_q^{2\vert1})$
the approach which Sudbery invented for $A_q(2)$ \cite{18} and also \cite{10}. Let us begin with the definition of the duality \cite{1}.

\begin{definition} \label{def5.1}
Two ${\mathbb Z}_2$-graded Hopf algebras ${\cal U}$ and ${\cal A}$ are said to be in duality if there exists a doubly non-degenerate bilinear form
$\langle,\rangle:{\cal U}\times{\cal A}\longrightarrow {\mathbb C}$, $(u,a)\mapsto \langle u,a\rangle$,
such that
\begin{subequations}\label{5.9}
\begin{align}
\langle uv, a\rangle &= \langle u \otimes v, \Delta_{\cal A}(a)\rangle, \quad \langle u, ab\rangle = \langle \Delta_{\cal U}(u), a \otimes b\rangle, \\
\langle u, {\bf 1}_{\cal A}\rangle &= \epsilon_{\cal U}(u), \quad \langle{\bf 1}_{\cal U}, a\rangle = \epsilon_{\cal A}(a), \\
\langle S_{\cal U}(u), a\rangle &= \langle u, S_{\cal A}(a)\rangle
\end{align}
\end{subequations}
for all $u, v \in {\cal U}$ and $a, b \in {\cal A}$.
\end{definition}

We note that equalities in (\ref{5.9}a) the pairing of ${\cal U}$ and ${\cal A}$ must be extended to one of ${\cal U} \otimes {\cal U}$ and
${\cal A} \otimes {\cal A}$ by setting
\begin{equation} \label{5.10}
\langle u \otimes v, a \otimes b\rangle = (-1)^{p(v)p(a)} \langle u,a\rangle\langle v,b\rangle.
\end{equation}
The pairing for elements of the ${\mathbb Z}_2$-graded Hopf algebras ${\cal U}$ and ${\cal A}$ follows from equations (\ref{5.9}) and the bilinear form inherited by the tensor product.

As a ${\mathbb Z}_2$-graded Hopf algebra ${\cal A}$ is generated by the elements $x$, $y$, $\theta$ and a basis is given by all monomials of the form $f = x^m y^n \theta^k$
where $m,n\in{\mathbb N}_0$ and $k\in\{0,1\}$. Let us denote the dual ${\mathbb Z}_2$-graded algebra by ${\cal U}_q$ and its generating elements by $X$, $Y$ and $\nabla$.

\begin{theorem} \label{thm5.3}
The commutation relations between the generators of the ${\mathbb Z}_2$-graded algebra ${\cal U}_q$ dual to ${\cal A}$ as follows
\begin{equation} \label{5.11}
[X,Y] = 0, \quad[X,\nabla] = \nabla, \quad [Y,\nabla] = -\nabla, \quad \nabla^2=0.
\end{equation}
\end{theorem}

\begin{proof}
The pairing is defined through the tangent vectors as follows
\begin{equation} \label{5.12}
\langle X, f\rangle = m \delta_{k,0}, \quad \langle Y, f\rangle = n\delta_{k,0}, \quad \langle \nabla, f\rangle  = \delta_{k,1}.
\end{equation}
We also have $\langle{\bf 1}_{\cal U}, f\rangle = \epsilon_{\cal A}(f) = \delta_{k,0}$. Using the defining relations (\ref{5.12}) one gets
\begin{equation*}
\langle X\nabla, f\rangle = m \delta_{k,1} \quad \mbox{and} \quad \langle\nabla X, f\rangle = (m-1) \delta_{k,1}
\end{equation*}
where differentiation is from the right as this is most suitable for differentiation in this basis. Thus, one obtains
\begin{equation*}
\langle X\nabla-\nabla X - \nabla, f\rangle = 0.
\end{equation*}
The other relations can be obtained similarly.
\end{proof}

\begin{theorem} \label{thm5.4}
The ${\mathbb Z}_2$-graded Hopf algebra structure of the ${\mathbb Z}_2$-graded algebra ${\cal U}_q$ is as follows:
\begin{align} \label{5.13}
\Delta_{\cal U}(X) &= X \otimes {\bf 1}_{\cal U} + {\bf 1}_{\cal U} \otimes X, & \epsilon_{\cal U}(X) &= 0, & S_{\cal U}(X) &= - X, \nonumber\\
\Delta_{\cal U}(Y) &= Y \otimes {\bf 1}_{\cal U} + {\bf 1}_{\cal U} \otimes Y, & \epsilon_{\cal U}(Y) &= 0, & S_{\cal U}(Y) &= - Y, \\
\Delta_{\cal U}(\nabla) &= \nabla \otimes {\bf 1}_{\cal U} + q^{X+Y}\otimes \nabla, & \epsilon_{\cal U}(\nabla) &= 0, & S_{\cal U}(\nabla) &= - q^{-(X+Y)}\nabla.\nonumber
\end{align}
\end{theorem}

\begin{proof}
We will only consider $\nabla$. The others can be shown similarly. So, we assume that $\Delta_{\cal U}(\nabla)= \nabla \otimes K_1 + K_2 \otimes \nabla$. Then
\begin{equation*}
\langle\Delta_{\cal U}(\nabla), x^my^n\otimes \theta\rangle
= \langle\nabla,x^my^n><K_1,\theta\rangle + \langle K_2,x^my^n\rangle\langle\nabla,\theta\rangle
= \langle K_2,x^my^n\rangle
\end{equation*}
and with (\ref{5.10})
\begin{equation*}
\langle\Delta_{\cal U}(\nabla), \theta\otimes x^my^n\rangle
= \langle\nabla,\theta\rangle\langle K_1,x^my^n\rangle - \langle K_2,\theta\rangle\langle\nabla,x^my^n\rangle \\
= \langle K_1,x^my^n\rangle,
\end{equation*}
so that we have $K_2=q^{X+Y}K_1$. To obtain the desired result we take $K_1={\bf 1}_{\cal U}$. The action of $\epsilon_{\cal U}$ and $S_{\cal U}$ on $\nabla$ follows from the identities $m\circ(\epsilon_{\cal U}\otimes{\rm id}_{\cal U})\circ\Delta_{\cal U}={\rm id}_{\cal U}$ and
$m\circ(S_{\cal U}\otimes{\rm id}_{\cal U})\circ\Delta_{\cal U} = \epsilon_{\cal U}$, respectively.
\end{proof}

Finally, let us write the relations (\ref{5.11}) as
\begin{equation} \label{5.14}
[\tilde{X},\tilde{Y}] = 0, \quad[\tilde{X},\nabla] = 0, \quad [\tilde{Y},\nabla] = 2\nabla, \quad \nabla^2=0.
\end{equation}
by choosing $\tilde{X}=X+Y$ and  $\tilde{Y}=X-Y$. We can now transform this algebra to the form obtained in Subsections 5.1 and 5.2 by having the following definitions:
\begin{equation*}
\tilde{X} \leftrightarrow N, \quad \tilde{Y} \leftrightarrow M, \quad \nabla \leftrightarrow \Phi q^{\tilde{X}}; \quad
\tilde{X} \leftrightarrow -(H+K), \quad \tilde{Y} \leftrightarrow H-K, \quad \nabla \leftrightarrow \Theta
\end{equation*}
which are consistent with the commutation relations and the Hopf structures.


\begin{thebibliography}{99}

\bibitem{1} E. Abe, {\em Hopf Agebras}, Cambridge Tracts in Math. N 74, Cambridge University Press, 1980.

\bibitem{2} T. Brzezinski and S. Majid,  {\em A class of bicovariant differential calculi on Hopf algebras}, Letters in Mathematical Physics, 26 (1992), 67--78.

\bibitem{3} S. Celik, {\em Differential geometry of the $q$-superplane}, Journal of Physics A: Mathematical and General, 31 (1998), 9695--9701.

\bibitem{4} S.A. Celik and S. Celik, {\em Differential geometry of the $q$-plane}, International Journal of Modern Physics A, 15 (2000), 3237--3243.

\bibitem{5} S. Celik, {\em Differential geometry of Z$_3$-graded quantum superplane}, Journal of Physics A: Mathematical and General, 35 (2002), 4257--4268.

\bibitem{6} S. Celik, {\em Cartan calculi on the quantum superplane} Journal of Mathematical Physics, 47 (2006), 083501:1--16.

\bibitem{7} S. Celik, {\em Bicovariant differential calculus on the quantum superspace ${\mathbb R}_q(1\vert2)$}, Journal of Algebras and its Applications, 15  (2016), 1650172:1--17.

\bibitem{8} S.A. Celik, {\em Differential calculi on super-Hopf algebra ${\mathbb F}({\mathbb R}_q(1\vert2))$}, Advanced in Applied Clifford Algebras, 28 (2018), 85:1--16.

\bibitem{9} S. Celik, {\em Covariant differential calculi on quantum symplectic superspace $SP_q^{1\vert2}$}, Journal of Mathematical Physics, 58 (2017), 023508: 1--15.

\bibitem{10} V.K. Dobrev, {\em Duality for the matrix quantum group ${\rm GL}_{p,q}(2,{\mathbb C})$}, Journal of Mathematical Physics, 33 (1992), 3419--3430.

\bibitem{11} H. Fakhri and S. Laheghi, {\em $GL_{r,s}(n)$-covariant differential calculi on the quantum $n$-space}, Advanced in Applied Clifford Algebras, 29 (2019), 52:1--10.

\bibitem{12} U. Hermisson, {\em Construction of covariant differential calculi on quantum homogeneous spaces}, Letters in Mathematical Physics, 46 (1998), 313--322.

\bibitem{13} T. Kobayashi and T. Uematsu, {\em Differential calculus on the quantum superspace and deformation of phase space}, Zeitschrift für Physik C: Particles and Fields, 56 (1992), 193--199.

\bibitem{14} Y.I. Manin, {\em Quantum Groups and Noncommutative Geometry}, Springer Nature Switzerland AG, 2018 (Second Edition).

\bibitem{15} Y.I. Manin, {\em Multiparametric quantum deformation of the general linear supergroup}, Communications in Mathematical Physics, 123 (1989), 163--175.

\bibitem{16} K. Schmuedgen, {\em On the Construction of covariant differential calculi on quantum homogeneous spaces}, Journal of Geometry and Physics, 30 (1999), 23--47.

\bibitem{17} S. Soni, {\em Differential calculus on the quantum superplane}, Journal of Physics A: Mathematical and General, 24 (1990), 619--624.

\bibitem{18} A. Sudbery, {\em Non-commuting coordinates and differential operators}, In: Curtrigh T et al (editors). Proceedings of the Argonne Workshop on Quantum Groups. Argone National Laboratory, (1990), 33--51.

\bibitem{19} J. Wess and B. Zumino, {\em Covariant differential calculus on the quantum hyperplane}, Nucleer Physics B, 18 (1990), 302--312.

\bibitem{20} S.L. Woronowicz, {\em Twisted $SU_q(2)$ group. An example of a non-commutative differential calculus}, Publications of the Research Institute for Mathematical Sciences, 23 (1987), 117--184.

\bibitem{21} S.L. Woronowicz, {\em Differential calculus on compact matric pseudogroups (quantum groups)}, Communications in Mathematical Physics, 122 (1989), 125--170.

\end{thebibliography}
\end{document}